\documentclass[reqno]{amsart}
\usepackage{amssymb,hyperref}

\newtheorem{thm}{Theorem}
\newtheorem{lem}{Lemma}

\renewcommand{\Re}{\mathrm{Re}}

\title
[Some sufficient problems for strongly close-to-convex of order $\mu$]
{Some sufficient problems for \\
strongly close-to-convex of order $\mu$}

\author{Hitoshi Shiraishi}
\address{Hitoshi Shiraishi \newline
Department of Mathematics \newline
Kinki University \newline
Higashi-Osaka, Osaka 577-8502, Japan}
\email{shiraishi@math.kindai.ac.jp}

\author{Shigeyoshi Owa}
\address{Shigeyoshi Owa \newline
Department of Mathematics \newline
Kinki University \newline
Higashi-Osaka, Osaka 577-8502, Japan}
\email{owa@math.kindai.ac.jp}


\date{}

\begin{document}

\begin{abstract}
For analytic functions $f(z)$ in the open unit disk $\mathbb{U}$ with $f(0)=f'(0)-1=0$,
a class $\mathcal{STC}(\mu)$ is defined.
The object of the present paper is to discuss some sufficient problems for $f(z)$ to be strongly close-to-convex of order $\mu$ in $\mathbb{U}$.
\end{abstract}

\begin{flushleft}
This paper was published in the journal: \\
Gen. Math. {\bf 17} (2009), No. 4, 157--169. \\
{\small \url{http://depmath.ulbsibiu.ro/genmath/gm/vol17nr4/13_Shiraishi/13_Shiraishi.pdf}}
\end{flushleft}
\hrule

\

\

\maketitle

\section{Introduction}

\

Let $\mathcal{A}_{n_j}$ denote the class of functions
$$
f(z)=z+a_{n_j+1}z^{n_j+1}+a_{n_j+2}z^{n_j+2}+ \ldots
\qquad(n=1,2,3,\ldots\ ;\ j=1,2)
$$
that are analytic in the open unit disk $\mathbb{U}=\{z \in \mathbb{C}:|z|<1\}$
and $\mathcal{A}=\mathcal{A}_1$.

We denote by $\mathcal{S}$ the subclass of $\mathcal{A}_n$ consisting of all univalent functions $f(z)$ in $\mathbb{U}$.

\

Let $\mathcal{S^{*}}(\alpha)$ be defined by
$$
\mathcal{S^{*}}(\alpha)=\left\{f(z)\in\mathcal{A}_n:\Re\left(\frac{zf'(z)}{f(z)}\right)>\alpha,\ 0 \leqq {}^{\exists}\alpha < 1\right\}.
$$
We denote by $\mathcal{S}^{*} = \mathcal{S}^{*}(0)$.

Also,
let $\mathcal{STC}(\mu)$ be defined by
$$
\mathcal{STC}(\mu)=\left\{f(z)\in\mathcal{A}_n:\Re\left(\left(\frac{zf'(z)}{g(z)}\right)^{\frac{1}{\mu}}\right)>0,\ 0 < {}^{\exists}\mu \leqq 1,\ {}^{\exists}g(z)\in\mathcal{S}^{*}\right\}.
$$
A function $f(z)\in\mathcal{STC}(\mu)$ is said to be strongly close-to-convex of order $\mu$ in $\mathbb{U}$.

\

The basic tool in proving our results is the following lemma due to Jack \cite{m1ref1}
(also, due to Miller and Mocanu \cite{m1ref2}).

\

\begin{lem} \label{jack} \quad
Let the function $w(z)$ defined by
$$
w(z)=a_nz^n+a_{n+1}z^{n+1}+a_{n+2}z^{n+2}+ \ldots
\qquad(n=1,2,3,\ldots)
$$
be analytic in $\mathbb{U}$ with $w(0)=0$.
If $\left|w(z)\right|$ attains its maximum value on the circle $\left|z\right|=r$ at a point $z_{0}\in\mathbb{U}$,
then there exists a real number $k \geqq n$ such that
$$
\frac{z_{0}w'(z_{0})}{w(z_{0})}=k.
$$
\end{lem}

\

\section{Main results}

\

Applying Lemma \ref{jack},
we drive the following results for $\mathcal{STC}(\mu)$.

\

\begin{thm} \label{d3thm1} \quad
If $f(z)\in\mathcal{A}_{n_1}$ satisfies
$$
\left| \left( \frac{zf'(z)}{g(z)} \right)^{\frac{1}{\mu}}-1\right|^{\beta}
\left( \Re\left( \delta +1 +\frac{zf''(z)}{f'(z)} -\frac{zg'(z)}{g(z)} \right) \right)^{\gamma}
< \left( \Re(\delta)+\frac{\mu n}{2} \right)^{\gamma}
\ (z\in\mathbb{U})
$$
for some real $\beta\geqq0$, $\gamma\geqq0$
such that $\beta + \gamma>0$,
some real $0< \mu \leqq 1$,
some complex $\delta$
with $\Re(\delta)>-\dfrac{\mu n}{2}$,
and for some $g(z)\in\mathcal{A}_{n_2}\cap\mathcal{S}^{*}$
where $n=\min\{n_1,n_2\}$,
then
$$
\left| \left(\frac{zf'(z)}{g(z)}\right)^\frac{1}{\mu} -1 \right|<1
\qquad(z\in\mathbb{U}).
$$

This means that $f(z)\in\mathcal{STC}(\mu)$.
\end{thm}

\

\begin{proof}\quad
Let us define $w(z)$ by
\begin{align}
w(z)
&= \left( \frac{zf'(z)}{g(z)} \right)^{\frac{1}{\mu}} -1
\qquad(z\in\mathbb{U})
\label{d3thm1eq1}\\
&= b_nz^n+b_{n+1}z^{n+1}+ \ldots \nonumber
\end{align}
where $n=\min\{n_1,n_2\}$.

Then, clearly, $w(z)$ is analytic in $\mathbb{U}$ and $w(0)=0$.
Differentiating both sides in (\ref{d3thm1eq1}),
we obtain
$$
1+\frac{zf''(z)}{f'(z)}-\frac{zg'(z)}{g(z)}=\frac{\mu zw'(z)}{w(z)+1},
$$
and therefore,
\begin{multline*}
\left| \left( \frac{zf'(z)}{g(z)} \right)^{\frac{1}{\mu}}-1\right|^{\beta}
\left( \Re\left( \delta +1 +\frac{zf''(z)}{f'(z)} -\frac{zg'(z)}{g(z)} \right) \right)^{\gamma} \\
= |w(z)|^{\beta} \left( \Re\left( \delta +\frac{\mu zw'(z)}{w(z)+1} \right) \right)^{\gamma}
< \left( \Re(\delta)+\frac{\mu n}{2} \right)^{\gamma}
\qquad(z\in\mathbb{U}).
\end{multline*}

If there exists a point $z_0 \in \mathbb{U}$ such that
$$
\max_{\left| z \right| \leqq \left| z_{0} \right|} \left| w(z) \right|
= \left| w(z_{0}) \right|
= 1,
$$
then Lemma \ref{jack} gives us that $w(z_0)=e^{i\theta}$ and $z_0 w'(z_0)=kw(z_0)$ ($k \geqq n$).

For such a point $z_0$,
we have
\begin{align*}
& \left| \left( \frac{z_0 f'(z_0)}{g(z_0)} \right)^{\frac{1}{\mu}}-1\right|^{\beta}
\left( \Re\left( \delta +1 +\frac{z_0 f''(z_0)}{f'(z_0)} -\frac{z_0 g'(z_0)}{g(z_0)} \right) \right)^{\gamma} \\
& = |w(z_0)|^{\beta} \left( \Re\left( \delta +\frac{\mu z_0 w'(z_0)}{w(z_0)+1} \right) \right)^{\gamma} \\
& = \left( \Re\left( \delta +\frac{\mu k w(z_0)}{w(z_0)+1} \right) \right)^{\gamma} \\
& = \left( \Re\left( \delta +\frac{\mu k}{2} \left( 1+i\tan\frac{\theta}{2}\right) \right) \right)^{\gamma} \\
& = \left( \Re(\delta) +\frac{\mu k}{2} \right)^{\gamma} \\
& \geqq \left( \Re(\delta) +\frac{\mu n}{2} \right)^{\gamma}.
\end{align*}

This contradicts our condition in the theorem.
Therefore, there is no $z_0 \in \mathbb{U}$ such that $|w(z_0)|=1$.
This means that $|w(z)|<1$ for all $z \in \mathbb{U}$.
It follows that
$$
\left| \left(\frac{zf'(z)}{g(z)}\right)^\frac{1}{\mu} -1 \right|<1
\qquad(z\in\mathbb{U})
$$
so that
$f(z) \in \mathcal{STC}(\mu)$.
\end{proof}

\

We also derive

\

\begin{thm} \label{d3thm2} \quad
If $f(z)\in\mathcal{A}_{n_1}$ satisfies
$$
\left| \left( \frac{zf'(z)}{g(z)} \right)^{\frac{1}{\mu}}-1\right|^{\beta}
\left( \Re \left( \delta +1 +\frac{zf''(z)}{f'(z)} \right) \right)^{\gamma}
\leqq \left( \Re(\delta) +\alpha +\frac{\mu n}{2} \right)^{\gamma}
\qquad(z\in\mathbb{U})
$$
for some real $\beta\geqq0$, $\gamma\geqq0$
such that $\beta + \gamma>0$,
some real $0< \mu \leqq 1$, $0 \leqq \alpha <1$,
some complex $\delta$
with $\Re(\delta)>-\dfrac{\mu n}{2}-\alpha$,
and for some $g(z)\in\mathcal{A}_{n_2}\cap\mathcal{S}^{*}(\alpha)$
where $n=\min\{n_1,n_2\}$,
then
$$
\left| \left(\frac{zf'(z)}{g(z)}\right)^\frac{1}{\mu} -1 \right|<1
\qquad(z\in\mathbb{U}),
$$
which shows that $f(z)\in\mathcal{STC}(\mu)$.
\end{thm}

\

\begin{proof}\quad
Define $w(z)$ in $\mathbb{U}$ by
\begin{align}
w(z)
&= \left( \frac{zf'(z)}{g(z)} \right)^{\frac{1}{\mu}} -1
\qquad(z\in\mathbb{U})
\label{d3thm2eq1}\\
&= b_nz^n+b_{n+1}z^{n+1}+ \ldots \nonumber
\end{align}
where $n=\min\{n_1,n_2\}$.

Evidently, $w(z)$ is analytic in $\mathbb{U}$ and $w(0)=0$.
Differentiating (\ref{d3thm2eq1}) logarithmically and simplyfing,
we have
$$
1+\frac{zf''(z)}{f'(z)}=\frac{zg'(z)}{g(z)}+\frac{\mu zw'(z)}{w(z)+1},
$$
and hence,
\begin{multline*}
\left| \left( \frac{zf'(z)}{g(z)} \right)^{\frac{1}{\mu}}-1\right|^{\beta}
\left( \Re\left( \delta +1 +\frac{zf''(z)}{f'(z)} \right) \right)^{\gamma} \\
= |w(z)|^{\beta} \left( \Re\left( \delta +\frac{zg'(z)}{g(z)} +\frac{\mu zw'(z)}{w(z)+1} \right) \right)^{\gamma}
\leqq \left( \Re(\delta) +\alpha +\frac{\mu n}{2} \right)^{\gamma}
\qquad(z\in\mathbb{U}).
\end{multline*}

If there exists a point $z_0 \in \mathbb{U}$ such that
$$
\max_{\left| z \right| \leqq \left| z_{0} \right|} \left| w(z) \right|
= \left| w(z_{0}) \right|
= 1,
$$
then Lemma \ref{jack} gives us that $w(z_0)=e^{i\theta}$ and $z_0 w'(z_0)=kw(z_0)$ ($k \geqq n$).

For such a point $z_0$,
we have
\begin{align*}
&\left| \left( \frac{z_0 f'(z_0)}{g(z_0)} \right)^{\frac{1}{\mu}}-1\right|^{\beta}
\left( \Re\left( \delta +1 +\frac{z_0 f''(z_0)}{f'(z_0)} \right) \right)^{\gamma} \\
& = |w(z_0)|^{\beta} \left( \Re\left( \delta +\frac{z_0 g'(z_0)}{g(z_0)} +\frac{\mu z_0 w'(z_0)}{w(z_0)+1} \right) \right)^{\gamma} \\
& = \left( \Re\left( \delta +\frac{z_0 g'(z_0)}{g(z_0)} +\frac{\mu k w(z_0)}{w(z_0)+1} \right) \right)^{\gamma} \\
& = \left( \Re\left( \delta +\frac{z_0 g'(z_0)}{g(z_0)} +\frac{\mu k}{2} \left( 1+i\tan\frac{\theta}{2}\right) \right) \right)^{\gamma} \\
& = \left( \Re(\delta) + \Re \left( \frac{z_0 g'(z_0)}{g(z_0)} \right) +\frac{\mu k}{2} \right)^{\gamma} \\
& > \left( \Re(\delta) +\alpha +\frac{\mu n}{2} \right)^{\gamma}.
\end{align*}

This contradicts our condition in the theorem.
Therefore, there is no $z_0 \in \mathbb{U}$ such that $|w(z_0)|=1$.
This means that $|w(z)|<1$ for all $z \in \mathbb{U}$.
This implies that
$$
\left| \left(\frac{zf'(z)}{g(z)}\right)^\frac{1}{\mu} -1 \right|<1
\qquad(z\in\mathbb{U})
$$
so that
$f(z) \in \mathcal{STC}(\mu)$.
\end{proof}

\

We consider a new appplication for Lemma \ref{jack}.
Our new application is follows.

\

\begin{thm} \label{d3thm3} \quad
If $f(z)\in\mathcal{A}_{n_1}$ satisfies
$$
\left| \left( \frac{zf'(z)}{g(z)} \right)^{\frac{1}{\mu}}-1\right|^{\beta}
\left| \delta +1 +\frac{zf''(z)}{f'(z)} -\frac{zg'(z)}{g(z)} \right|^{\gamma}
< \rho^{\beta} \left( \delta +\frac{\mu \rho n}{1+\rho} \right)^{\gamma}
\qquad(z\in\mathbb{U})
$$
for some real $\beta\geqq0$, $\gamma\geqq0$
such that $\beta + \gamma>0$,
some real $0< \mu \leqq 1$, $\delta>0$, $\rho$
with $\rho>\sqrt[]{\dfrac{\delta}{\delta+\mu n}}$,
and for some $g(z)\in\mathcal{A}_{n_2}\cap\mathcal{S}^{*}$
where $n=\min\{n_1,n_2\}$,
then
$$
\left| \left(\frac{zf'(z)}{g(z)}\right)^\frac{1}{\mu} -1 \right|<\rho
\qquad(z\in\mathbb{U}).
$$

In addition for $\rho<1$,
we have $f(z)\in\mathcal{STC}(\mu)$.
\end{thm}

\

\begin{proof}\quad
Defining the function $w(z)$ by
\begin{align*}
w(z)
&= \left( \frac{zf'(z)}{g(z)} \right)^{\frac{1}{\mu}} -1
\qquad(z\in\mathbb{U})\\
&= b_nz^n+b_{n+1}z^{n+1}+ \ldots
\end{align*}
where $n=\min\{n_1,n_2\}$,
we have that $w(z)$ is analytic in $\mathbb{U}$ with $w(0)=0$.
Since,
$$
1+\frac{zf''(z)}{f'(z)}-\frac{zg'(z)}{g(z)}=\frac{\mu zw'(z)}{w(z)+1},
$$
we obtain that
\begin{align*}
\left| \left( \frac{zf'(z)}{g(z)} \right)^{\frac{1}{\mu}}-1\right|^{\beta}
\left| \delta +1 +\frac{zf''(z)}{f'(z)} -\frac{zg'(z)}{g(z)} \right|^{\gamma}
& = |w(z)|^{\beta} \left| \delta +\frac{\mu zw'(z)}{w(z)+1} \right|^{\gamma} \\
& < \rho^{\beta} \left( \delta +\frac{\mu \rho n}{1+\rho} \right)^{\gamma}
\qquad(z\in\mathbb{U}).
\end{align*}

If there exists a point $z_0 \in \mathbb{U}$ such that
$$
\max_{\left| z \right| \leqq \left| z_{0} \right|} \left| w(z) \right|
= \left| w(z_{0}) \right|
= \rho,
$$
then Lemma \ref{jack} gives us that $w(z_0)= \rho e^{i\theta}$ and $z_0 w'(z_0)=kw(z_0)$ ($k \geqq n$).

Thus we have
\begin{align*}
& \left| \left( \frac{z_0 f'(z_0)}{g(z_0)} \right)^{\frac{1}{\mu}}-1\right|^{\beta}
\left| \delta +1 +\frac{z_0 f''(z_0)}{f'(z_0)} -\frac{z_0 g'(z_0)}{g(z_0)} \right|^{\gamma} \\
& = |w(z_0)|^{\beta} \left| \delta +\frac{\mu z_0 w'(z_0)}{w(z_0)+1} \right|^{\gamma} \\
& = \rho^\beta \left| \delta +\frac{\mu k w(z_0)}{w(z_0)+1} \right|^{\gamma} \\
& = \rho^\beta \left| \delta +\frac{\mu \rho k (\rho + \cos\theta)}{\rho^2 +1 +2\rho\cos\theta} +i\frac{\mu \rho k \sin\theta}{\rho^2 +1 +2\rho\cos\theta} \right|^{\gamma} \\
& = \rho^\beta \left( \delta^2 +\mu \delta k +\frac{\mu \delta k (\rho^2 -1) +\mu^2 \rho^2 k^2}{\rho^2 +1 +2\rho\cos\theta} \right)^{\frac{\gamma}{2}} \\
& \geqq \rho^\beta \left( \delta^2 +\mu \delta k +\frac{\mu \delta k (\rho^2 -1) +\mu^2 \rho^2 k^2}{\rho^2 +1 +2\rho} \right)^{\frac{\gamma}{2}} \\
& = \rho^{\beta} \left( \delta +\frac{\mu \rho k}{1+\rho} \right)^{\gamma} \\
& \geqq \rho^{\beta} \left( \delta +\frac{\mu \rho n}{1+\rho} \right)^{\gamma}.
\end{align*}

This contradicts our condition in the theorem.
Therefore, there is no $z_0 \in \mathbb{U}$ such that $|w(z_0)|=\rho$.
This means that $|w(z)|<\rho$ for all $z \in \mathbb{U}$.
\end{proof}

\

We also consider a new aplocstion for Lemma \ref{jack}.

\

\begin{thm} \label{d3thm4} \quad
If $f(z)\in\mathcal{A}_{n_1}$ satisfies
$$
\left| \left( \frac{g(z)}{zf'(z)} \right)^{\frac{1}{\mu}}-1\right|^{\beta}
\left( \Re\left( \delta -1 -\frac{zf''(z)}{f'(z)} +\frac{zg'(z)}{g(z)} \right) \right)^{\gamma}
< \left( \Re(\delta)+\frac{\mu n}{2} \right)^{\gamma}
\ (z\in\mathbb{U})
$$
for some real $\beta\geqq0$, $\gamma\geqq0$
such that $\beta + \gamma>0$,
some real $0< \mu \leqq 1$,
some complex $\delta$
with $\Re(\delta)>-\dfrac{\mu n}{2}$,
and for some $g(z)\in\mathcal{A}_{n_2}\cap\mathcal{S}^{*}$
where $n=\min\{n_1,n_2\}$,
then
$$
\left| \left( \frac{g(z)}{zf'(z)} \right)^\frac{1}{\mu} -1 \right|<1
\qquad(z\in\mathbb{U}).
$$

This means that $f(z)\in\mathcal{STC}(\mu)$.
\end{thm}

\

\begin{proof}\quad
Let us define the function $w(z)$ by
\begin{align}
w(z)
&= \left( \frac{g(z)}{zf'(z)} \right)^{\frac{1}{\mu}} -1
\qquad(z\in\mathbb{U})\\
&= b_nz^n+b_{n+1}z^{n+1}+ \ldots \nonumber
\end{align}
where $n=\min\{n_1,n_2\}$.

Clearly, $w(z)$ is analytic in $\mathbb{U}$ with $w(0)=0$.
We want to prove that $|w(z)|<1$ in $\mathbb{U}$.
Since,
$$
-1-\frac{zf''(z)}{f'(z)}+\frac{zg'(z)}{g(z)}=\frac{\mu zw'(z)}{w(z)+1},
$$
we see that
\begin{multline*}
\left| \left( \frac{g(z)}{zf'(z)} \right)^{\frac{1}{\mu}}-1\right|^{\beta}
\left( \Re\left( \delta -1 -\frac{zf''(z)}{f'(z)} +\frac{zg'(z)}{g(z)} \right) \right)^{\gamma} \\
= |w(z)|^{\beta} \left( \Re\left( \delta +\frac{\mu zw'(z)}{w(z)+1} \right) \right)^{\gamma}
< \left( \Re(\delta)+\frac{\mu n}{2} \right)^{\gamma}
\qquad(z\in\mathbb{U}).
\end{multline*}

If there exists a point $z_0 \in \mathbb{U}$ such that
$$
\max_{\left| z \right| \leqq \left| z_{0} \right|} \left| w(z) \right|
= \left| w(z_{0}) \right|
= 1,
$$
then Lemma \ref{jack} gives us that $w(z_0)=e^{i\theta}$ and $z_0 w'(z_0)=kw(z_0)$ ($k \geqq n$).

Thus we have
\begin{align*}
& \left| \left( \frac{g(z_0)}{z_0 f'(z_0)} \right)^{\frac{1}{\mu}}-1\right|^{\beta}
\left( \Re\left( \delta -1 -\frac{z_0 f''(z_0)}{f'(z_0)} +\frac{z_0 g'(z_0)}{g(z_0)} \right) \right)^{\gamma} \\
& = |w(z_0)|^{\beta} \left( \Re\left( \delta +\frac{\mu z_0 w'(z_0)}{w(z_0)+1} \right) \right)^{\gamma} \\
& = \left( \Re\left( \delta +\frac{\mu k w(z_0)}{w(z_0)+1} \right) \right)^{\gamma} \\
& = \left( \Re\left( \delta +\frac{\mu k}{2} \left( 1+i\tan\frac{\theta}{2}\right) \right) \right)^{\gamma} \\
& = \left( \Re(\delta) +\frac{\mu k}{2} \right)^{\gamma} \\
& \geqq \left( \Re(\delta) +\frac{\mu n}{2} \right)^{\gamma}.
\end{align*}

This contradicts the condition in the theorem.
Therefore, there is no $z_0 \in \mathbb{U}$ such that $|w(z_0)|=1$.
This means that $|w(z)|<1$ for all $z \in \mathbb{U}$.
We conclude that
$$
\left| \left(\frac{g(z)}{zf'(z)}\right)^\frac{1}{\mu} -1 \right|<1
\qquad(z\in\mathbb{U})
$$
which implies that
$f(z) \in \mathcal{STC}(\mu)$.
\end{proof}

\

Finally,
we derive

\

\begin{thm} \label{d3thm5} \quad
If $f(z)\in\mathcal{A}_{n_1}$ satisfies
$$
\left| \left( \frac{g(z)}{zf'(z)} \right)^{\frac{1}{\mu}}-1\right|^{\beta}
\left| \delta -1 -\frac{zf''(z)}{f'(z)} +\frac{zg'(z)}{g(z)} \right|^{\gamma}
< \rho^{\beta} \left( \delta +\frac{\mu \rho n}{1+\rho} \right)^{\gamma}
\qquad(z\in\mathbb{U})
$$
for some real $\beta\geqq0$, $\gamma\geqq0$
such that $\beta + \gamma>0$,
some real $0< \mu \leqq 1$, $\delta>0$, $\rho$
with $\rho>\sqrt[]{\dfrac{\delta}{\delta+\mu n}}$,
and for some $g(z)\in\mathcal{A}_{n_2}\cap\mathcal{S}^{*}$
where $n=\min\{n_1,n_2\}$,
then
$$
\left| \left(\frac{g(z)}{zf'(z)}\right)^\frac{1}{\mu} -1 \right|<\rho
\qquad(z\in\mathbb{U}).
$$

In addition,
if $\rho<1$,
then we have $f(z)\in\mathcal{STC}(\mu).$
\end{thm}

\

\begin{proof}\quad
Let us define $w(z)$ by
\begin{align}
w(z)
&= \left( \frac{g(z)}{zf'(z)} \right)^{\frac{1}{\mu}} -1
\qquad(z\in\mathbb{U})
\label{d3thm5eq1}\\
&= b_nz^n+b_{n+1}z^{n+1}+ \ldots \nonumber
\end{align}
where $n=\min\{n_1,n_2\}$.

Then, we have that $w(z)$ is analytic in $\mathbb{U}$ and $w(0)=0$.
Differenciating (\ref{d3thm5eq1}) in both side logarithmically and simplifying,
we obtain
$$
-1-\frac{zf''(z)}{f'(z)}+\frac{zg'(z)}{g(z)}=\frac{\mu zw'(z)}{w(z)+1},
$$
and hence,
\begin{align*}
\left| \left( \frac{g(z)}{zf'(z)} \right)^{\frac{1}{\mu}}-1\right|^{\beta}
\left| \delta -1 -\frac{zf''(z)}{f'(z)} +\frac{zg'(z)}{g(z)} \right|^{\gamma}
& = |w(z)|^{\beta} \left| \delta +\frac{\mu zw'(z)}{w(z)+1} \right|^{\gamma} \\
& < \rho^{\beta} \left( \delta +\frac{\mu \rho n}{1+\rho} \right)^{\gamma}
\qquad(z\in\mathbb{U}).
\end{align*}

If there exists a point $z_0 \in \mathbb{U}$ such that
$$
\max_{\left| z \right| \leqq \left| z_{0} \right|} \left| w(z) \right|
= \left| w(z_{0}) \right|
= \rho,
$$
then Lemma \ref{jack} gives us that $w(z_0)= \rho e^{i\theta}$ and $z_0 w'(z_0)=kw(z_0)$ ($k \geqq n$).

Therefore,
we have
\begin{align*}
& \left| \left( \frac{g(z_0)}{z_0 f'(z_0)} \right)^{\frac{1}{\mu}}-1\right|^{\beta}
\left| \delta -1 -\frac{z_0 f''(z_0)}{f'(z_0)} +\frac{z_0 g'(z_0)}{g(z_0)} \right|^{\gamma} \\
& = |w(z_0)|^{\beta} \left| \delta +\frac{\mu z_0 w'(z_0)}{w(z_0)+1} \right|^{\gamma} \\
& = \rho^\beta \left| \delta +\frac{\mu k w(z_0)}{w(z_0)+1} \right|^{\gamma} \\
& = \rho^\beta \left| \delta +\frac{\mu \rho k (\rho + \cos\theta)}{\rho^2 +1 +2\rho\cos\theta} +i\frac{\mu \rho k \sin\theta}{\rho^2 +1 +2\rho\cos\theta} \right|^{\gamma} \\
& = \rho^\beta \left( \delta^2 +\mu \delta k +\frac{\mu \delta k (\rho^2 -1) +\mu^2 \rho^2 k^2}{\rho^2 +1 +2\rho\cos\theta} \right)^{\frac{\gamma}{2}} \\
& \geqq \rho^\beta \left( \delta^2 +\mu \delta k +\frac{\mu \delta k (\rho^2 -1) +\mu^2 \rho^2 k^2}{\rho^2 +1 +2\rho} \right)^{\frac{\gamma}{2}} \\
& = \rho^{\beta} \left( \delta +\frac{\mu \rho k}{1+\rho} \right)^{\gamma} \\
& \geqq \rho^{\beta} \left( \delta +\frac{\mu \rho n}{1+\rho} \right)^{\gamma}.
\end{align*}

This contradicts our condition in the theorem.
Therefore, there is no $z_0 \in \mathbb{U}$ such that $|w(z_0)|=\rho$.
This means that $|w(z)|<\rho$ for all $z \in \mathbb{U}$.
\end{proof}

\


\begin{thebibliography}{}

\bibitem{m1ref1}
I. S. Jack,
{\it Functions starlike and convex of order $\alpha$},
J. London Math. Soc. {\bf 3}(1971), 469--474.

\bibitem{m1ref2}
S. S. Miller and P. T. Mocanu,
{\it Second-order differential inequalities in the complex plane},
J. Math. Anal. Appl. {\bf 65}(1978), 289--305.

\bibitem{m1ref3}
R. Singh and S. Singh,
{\it Some sufficient conditions for univalence and starlikeness},
Coll. Math. {\bf 47}(1982), 309--314.

\bibitem{m1ref0}
H. Shiraishi and S. Owa,
{\it Some sufficient problems for certain univalent functions},
Far East J. Math. Sci. {\bf 30}(2008), 147--155.

\end{thebibliography}
\end{document}